\documentclass[12pt]{article}                                               

\input{amssym.def}                                                           
\input{amssym.tex}                                                           
                                                                             
\begin{document}                                                             
\title{On a symmetry of complex and  real   multiplication}

\author{Igor ~V. ~Nikolaev
\footnote{Partially supported 
by NSERC.}}


\date{}
 \maketitle


\newtheorem{thm}{Theorem}
\newtheorem{lem}{Lemma}
\newtheorem{dfn}{Definition}
\newtheorem{rmk}{Remark}
\newtheorem{cor}{Corollary}
\newtheorem{cnj}{Conjecture}
\newtheorem{exm}{Example}



\begin{abstract}
It is proved that each lattice with  complex multiplication by $f\sqrt{-D}$ 
corresponds to  a  pseudo-lattice  with  real multiplication by $f'\sqrt{D}$,
where $f'$ is an integer defined by $f$.

\vspace{7mm}


\vspace{5mm}
\end{abstract}

\section{Introduction}
The paper continues a study of the  duality  between elliptic curves with
complex multiplication and noncommutative tori with real multiplication
initiated in \cite{Nik1};  let us introduce some notation and basic facts.    
Fix an irrational   number  $0<\theta <1$;   a 
{\it noncommutative torus} is the universal $C^*$-algebra $A_{\theta}$ generated by the  unitaries 
$u$ and $v$ satisfying the commutation  relation $vu=e^{2\pi i\theta}uv$ (Rieffel,  1981  \cite{Rie1}). 
Two such tori  are stably isomorphic (Morita equivalent)  whenever $A_{\theta}\otimes {\cal K}
\cong A_{\theta'}\otimes {\cal  K}$,  where ${\cal K}$ is the $C^*$-algebra of 
compact operators;  the isomorphism occurs  if and only if  $\theta'= (a\theta +b) / (c\theta +d)$,  
where  $a,b,c,d\in {\Bbb Z}$ and  $ad-bc=1$.
The K-theory of $A_{\theta}$ is  Bott periodic with  $K_0(A_{\theta})=K_1(A_{\theta})\cong {\Bbb Z}^2$.
The range of the trace on  projections of $A_{\theta}\otimes {\cal K}$ is a subset
$\Lambda={\Bbb Z}+{\Bbb Z}\theta$  of the real line (Rieffel,   1981  \cite{Rie1});   $\Lambda$ is called a 
pseudo-lattice (Manin,  2004   \cite{Man1}).   The torus $A_{\theta}$ is said to have  {\it real multiplication} 
if $\theta$ is a quadratic irrationality;  we shall denote the set of such algebras   by ${\cal A}_{RM}$.
The real multiplication entails existence of the non-trivial  endomorphisms  of $\Lambda$ coming from multiplication
by the real numbers -- hence the name.   If $D>1$ is a square-free  integer,
we shall write $A_{RM}^{(D,f)}$ to denote 
real multiplication by an order $R_f$  of conductor $f\ge 1$
in the  field ${\Bbb Q}(\sqrt{D})$;  each torus in ${\cal A}_{RM}$ can be written
in this form (Manin,  2004   \cite{Man1}).

Let ${\Bbb H}=\{x+iy\in {\Bbb C}~|~y>0\}$ be the
upper half-plane and for $\tau\in {\Bbb H}$ let ${\Bbb C}/({\Bbb Z}+{\Bbb Z}\tau)$
be a complex torus;  we routinely identify the latter with a non-singular elliptic curve via the
Weierstrass $\wp$ function (Silverman, 1994   [7, pp. 6-7]).  Recall that  two complex tori are 
isomorphic, whenever  $\tau'= (a\tau +b) / (c\tau +d)$,  where  $a,b,c,d\in {\Bbb Z}$ and 
$ad-bc=1$. If $\tau$ is an imaginary quadratic number,  the elliptic curve is said to have 
{\it complex multiplication};   in this case the  lattice $L={\Bbb Z}+{\Bbb Z}\tau$
admits non-trivial endomorphisms given as multiplication of $L$ by certain complex (quadratic) numbers.  
Elliptic curves with complex multiplication  are fundamental
and have long history in  number theory;  we shall denote the set of such curves  by ${\cal E}_{CM}$.
We  write $E_{CM}^{(-D,f)}$ to denote  the  elliptic curve with
complex multiplication by an order ${\goth R}_f$  of conductor $f\ge 1$ in the
imaginary quadratic  field ${\Bbb Q}(\sqrt{-D})$;  each curve in ${\cal E}_{CM}$ is isomorphic
to $E_{CM}^{(-D,f)}$ for some integers $D$ and $f$  (Silverman,  1994  [7,  pp. 95-96]).

There exists  a covariant functor  between  elliptic curves and  noncommutative tori;   
the functor maps  isomorphic curves to the stably isomorphic 
tori.  To give an idea,  let $\phi$ be a closed form on a topological torus;
the trajectories of $\phi$  define a measured foliation on the torus.  
By the Hubbard-Masur  theorem, such a foliation 
corresponds to a point $\tau\in {\Bbb H}$.  The map $F: {\Bbb H}\to\partial {\Bbb H}$
is defined by the formula $\tau\mapsto\theta=\int_{\gamma_2}\phi/\int_{\gamma_1}\phi$,
where $\gamma_1$ and $\gamma_2$ are generators of the first homology of the
torus.    The following is true: (i) ${\Bbb H}=\partial {\Bbb H}\times (0,\infty)$
is a trivial fiber bundle, whose projection map coincides with $F$;
(ii) $F$ is a functor, which maps isomorphic complex tori to
the stably isomorphic noncommutative tori.  We shall refer to $F$
as the {\it Teichm\"uller functor}.   It was proved in \cite{Nik1} that 
$F({\cal E}_{CM})\subseteq {\cal A}_{RM}$, i.e.   $F$  sends elliptic curves with complex 
multiplication to the noncommutative tori with real multiplication. 
Namely,  $F(E_{CM}^{(-D,f)})=A_{RM}^{(D,f')}$,  where $f'$ is the least integer satisfying 
equation $|Cl~(R_{f'})|=|Cl~({\goth R}_f)|$ for the class numbers of orders $R_{f'}$ and ${\goth R}_f$,
respectively;  the latter constraint is a necessary and sufficient condition for $A_{RM}^{(D,f')}$
to discern non-isomorphic curves  $E_{CM}^{(-D,f)}$ having the same endomorphism ring 
$R_f$.

Denote by $\Lambda_{RM}^{(D,f)}$  a pseudo-lattice corresponding to the torus $A_{RM}^{(D,f)}$;
the  $\Lambda_{RM}^{(D,f)}$ can be identified with  points of  the boundary $\partial {\Bbb H}$  of the
half-plane ${\Bbb H}$. 
Let $x,\bar x\in \Lambda_{RM}^{(D,f)}$ be a pair of the  conjugate quadratic irrationalities and 
consider  a geodesic half-circle through $x$ and $\bar x$:
\begin{equation}\label{eq0}
\widetilde \gamma (x,\bar x)={xe^{t\over 2}+i\bar x e^{-{t\over 2}}\over e^{t\over 2}+ie^{-{t\over 2}}},
\qquad -\infty\le t\le \infty.
\end{equation}
A Riemann surface $X$ is said to be {\it associated to}  $A_{RM}^{(D,f)}$,
if the covering of the geodesic spectrum of $X$  contains 
the set $\{\widetilde\gamma (x,\bar x) ~: ~\forall x\in \Lambda_{RM}^{(D,f)}\}$,
see definition \ref{dfn1};  such a surface will be denoted 
by  $X(A_{RM}^{(D,f)})$. 
 Our main result can be expressed as follows.   
\begin{thm}\label{thm1}
For every square-free  integer $D>1$ and integer $f\ge 1$  there exists a holomorphic 
map $F^{-1}: X(A_{RM}^{(D,f')})\to E_{CM}^{(-D,f)}$,  where
$F(E_{CM}^{(-D,f)})=A_{RM}^{(D,f')}$.
\end{thm}
The note is organized as follows. 
Section 2 is reserved for notation and preliminary facts.  
Theorem \ref{thm1} is proved in Section 3.

\section{Riemann surface $X(A_{RM}^{(D,f)})$}
Let $X$ be a Riemann surface;  consider the geodesic spectrum of $X$,
i.e. the set  $Spec~X$ consisting of all closed geodesics of  $X$.  
Recall that for the covering map ${\Bbb H}\to X$   each geodesic 
$\gamma\in Spec~X$  is the image of a geodesic  half-circle 
$\widetilde\gamma(x, x')\in {\Bbb H}$ with the endpoints  
$x\ne x'$.   Denote by  $\widetilde{Spec}~X\subset {\Bbb H}$  the set of 
geodesic half-circles covering the geodesic spectrum of  $X$. 
\begin{dfn}\label{dfn1}
We shall say that the Riemann surface $X$ is associated to the noncommutative
torus  $A_{RM}^{(D,f)}$, if $\{\widetilde\gamma (x,\bar x) ~: ~\forall x\in \Lambda_{RM}^{(D,f)}\}
\subset \widetilde{Spec}~X$;  the associated 
Riemann surface will be denoted by  $X(A_{RM}^{(D,f)})$.
 \end{dfn}
Let $N\ge 1$ be an integer;  by $\Gamma_1(N)$ we understand a subgroup 
of the modular group $SL_2({\Bbb Z})$ consisting of matrices of the form
\begin{equation}
\left\{\left(\matrix{a & b\cr c & d}\right)\in SL_2({\Bbb Z}) ~| ~ a,d\equiv~1~mod~N, 
~c\equiv~0~mod~N\right\};
\end{equation}
the corresponding Riemann surface ${\Bbb H}/\Gamma_1(N)$ will be denoted by $X_1(N)$. 
The following lemma links  $X(A_{RM}^{(D,f)})$ to $X_1(N)$.   
\begin{lem}\label{lm1}
$X(A_{RM}^{(D,f)})\cong X_1(fD)$.
\end{lem}
{\it Proof.}
Let $\Lambda_{RM}^{(D,f)}$ be a pseudo-lattice with  real 
multiplication by an order $R$ in the real quadratic number field
${\Bbb Q}(\sqrt{D})$; it is known, that $\Lambda_{RM}^{(D,f)}\subseteq R$
and $R={\Bbb Z}+(f\omega){\Bbb Z}$, where $f\ge 1$ is the conductor of $R$ and  
\begin{equation}\label{eq1}
\omega=\cases{{1+\sqrt{D}\over 2} & if $D\equiv 1 ~mod~4$,\cr
               \sqrt{D} & if $D\equiv 2,3 ~mod~4$,}
\end{equation}
see e.g. (Borevich \&  Shafarevich,  1988   [1,  pp. 130-131])
Recall that matrix $(a,b,c,d)\in SL_2({\Bbb Z})$ has a pair
of real fixed points $x$ and $\bar x$ if and only if $|a+d|>2$ (the hyperbolic matrix);
the fixed points can be found from the equation $x=(ax+b)(cx+d)^{-1}$  by the
formulas:
\begin{equation}\label{eq2}
x={a-d\over 2c}+\sqrt{{(a+d)^2-4\over 4c^2}}, \qquad
\bar x={a-d\over 2c}-\sqrt{{(a+d)^2-4\over 4c^2}}. 
\end{equation}

\bigskip
{\sf Case I.} If $D\equiv 1~mod~4$, then formula (\ref{eq1}) implies that
$R=(1+{f\over 2}){\Bbb Z}+{\sqrt{f^2D}\over 2}{\Bbb Z}$. If $x\in\Lambda_{RM}^{(D,f)}$
is fixed point of a transformation $(a,b,c,d)\in SL_2({\Bbb Z})$, then
formula (\ref{eq2}) implies:
\begin{equation}
\left\{
\begin{array}{ccc}
{a-d\over 2c} &=& (1+{f\over 2})z_1\\
{(a+d)^2-4\over 4c^2}  &=& {f^2D\over 4}z_2^2
\end{array}
\right.
\end{equation}
for some integer numbers $z_1$ and $z_2$.  
The second equation can be written in the form $(a+d)^2-4=c^2f^2Dz_2^2$;
we have therefore $(a+d)^2\equiv 4~mod~(fD)$ and $a+d\equiv\pm 2~mod~(fD)$.
Without loss of generality we assume $a+d\equiv 2~mod~(fD)$ since
matrix $(a,b,c,d)\in SL_2({\Bbb Z})$ can be multiplied by $-1$.  Notice that
the last equation admits a solution $a=d\equiv 1~mod~(fD)$.

The first equation yields us ${a-d\over c}=(2+f)z_1$, where $c\ne0$ since
the matrix $(a,b,c,d)$ is hyperbolic.  Notice that  $a-d\equiv 0~mod~(fD)$;
since the ratio ${a-d\over c}$ must be integer, we conclude that $c\equiv 0~mod~(fD)$.  
All together,  we get: 
\begin{equation}\label{eq4}
a\equiv 1~mod~(fD), \quad d\equiv  1~mod~(fD), \quad  c\equiv 0~mod~(fD).
\end{equation}

\bigskip
{\sf Case II.} If $D\equiv 2$ or $3~mod~4$, then
$R={\Bbb Z}+(\sqrt{f^2D})~{\Bbb Z}$. If $x\in\Lambda_{RM}^{(D,f)}$
is  fixed point of a  transformation $(a,b,c,d)\in SL_2({\Bbb Z})$, then
formula (\ref{eq2}) implies:
\begin{equation}
\left\{
\begin{array}{ccc}
{a-d\over 2c} &=& z_1\\
{(a+d)^2-4\over 4c^2}  &=& f^2Dz_2^2
\end{array}
\right.
\end{equation}
for some integer numbers $z_1$ and $z_2$.   The second equation
gives  $(a+d)^2-4=4c^2f^2Dz_2^2$;  therefore  $(a+d)^2\equiv 4~mod~(fD)$ and $a+d\equiv\pm 2~mod~(fD)$.
Again without loss of generality we assume $a+d\equiv 2~mod~(fD)$ since
matrix $(a,b,c,d)\in SL_2({\Bbb Z})$ can be  multiplied by $-1$.  
The last equation admits a solution $a=d\equiv 1~mod~(fD)$.

The first equation is ${a-d\over c}=2z_1$, where $c\ne0$.
Since  $a-d\equiv 0~mod~(fD)$ and the ratio ${a-d\over c}$ must be integer,
one concludes  that $c\equiv 0~mod~(fD)$.
 All together,  we get equations (\ref{eq4}).  
Since all possible cases are exhausted,  lemma \ref{lm1} follows.
$\square$
\begin{rmk}\label{rmk2}
There exist other finite index subgroups of $SL_2({\Bbb Z})$ whose
geodesic spectrum contains the set 
 $\{\widetilde\gamma (x,\bar x) ~: ~\forall x\in \Lambda_{RM}^{(D,f)}\}$;
however $\Gamma_1(fD)$ is a unique group with such a property 
among subgroups of the principal congruence group.
\end{rmk}
\begin{rmk}\label{rmk1}
Not all geodesics of $X_1(fD)$  have  form (\ref{eq0});  thus the set \linebreak
 $\{\widetilde\gamma (x,\bar x) ~: ~\forall x\in \Lambda_{RM}^{(D,f)}\}$ is strictly
included in the geodesic spectrum of modular curve $X_1(fD)$.
\end{rmk}

\section{Proof of theorem \ref{thm1}}
Recall, that $\Gamma(N):=\{(a,b,c,d)\in SL_2({\Bbb Z})~|~a,d\equiv~1~mod~N,~b,c\equiv~0~mod~N\}$
is called a {\it principal congruence group} of level $N$;  the corresponding (compact) Riemann surface will be
denoted by $X(N)={\Bbb H}/\Gamma(N)$. 
\begin{lem}\label{lm2}
 {\bf (Hecke)}
There exists a holomorphic map $X(fD)\to E_{CM}^{(-D,f)}$. 
\end{lem}
{\it Proof.} A detailed proof of this beautiful fact is given in
(Hecke,  1928   \cite{Hec1}).

To give an idea of the proof, let  $\goth R$ be an order of conductor $f\ge 1$
in the imaginary quadratic number field ${\Bbb Q}(\sqrt{-D})$; consider
an $L$-function attached to ${\goth R}$:
\begin{equation}
L(s, \psi)=\prod_{{\goth P}\subset {\goth R}}{1\over 1-{\psi({\goth P})\over N({\goth P})^s}},
\quad s\in {\Bbb C}, 
\end{equation}
where ${\goth P}$ is a prime ideal in ${\goth R}$, $N({\goth P})$ its norm
and $\psi$ a  Gr\"ossencharacter. 
A crucial observation (\S 1) says that the series $L(s, \psi)$ converges to a cusp form $w(s)$
of the principal congruence group $\Gamma(fD)$.

By the Deuring Theorem,   $L(E_{CM}^{(-D,f)},s)=L(s,\psi)L(s, \bar\psi)$, where $L(E_{CM}^{(-D,f)},s)$ is the
Hasse-Weil $L$-function of the elliptic curve and $\bar\psi$ a conjugate of the 
Gr\"ossencharacter,  see  (Silverman,   1994   [7,  p. 175]);
moreover $L(E_{CM}^{(-D,f)},s)=L(w, s)$,  where $L(w, s):=\sum_{n=1}^{\infty}{c_n\over n^s}$ and $c_n$
the Fourier coefficients of the cusp form $w(s)$.  In other words,   
$E_{CM}^{(-D,f)}$ is a modular elliptic curve.

One can now apply the modularity principle: if  $A_w$ is an abelian variety given by 
the periods of holomorphic differential $w(s)ds$ (and its conjugates)  on  
$X(fD)$, then the following diagram commutes

\bigskip   
\begin{picture}(300,110)(-80,-5)
\put(130,70){\vector(0,-1){35}}
\put(55,70){\vector(2,-1){53}}
\put(55,83){\vector(1,0){53}}
\put(128,20){$E_{CM}^{(-D,f)}$}
\put(10,80){$X(fD)$}
\put(125,80){$A_w$}
\put(80,90){$\iota$}
\put(140,55){$\pi$}
\end{picture}

\noindent
The holomorphic map  $X(fD)\to E_{CM}^{(-D,f)}$ is obtained as  a composition of the canonical
embedding $\iota: X(fD)\to A_w$ with the subsequent holomorphic projection 
$\pi:  A_w\to E_{CM}^{(-D,f)}$.
$\square$

\begin{lem}\label{lm3}
The functor $F$ acts by the formula
$E_{CM}^{(-D,f)}\mapsto A_{RM}^{(D,f)}$. 
\end{lem}
{\it Proof.} Let $L_{CM}$ be a lattice with complex multiplication by an order ${\goth R}={\Bbb Z}+(f\omega){\Bbb Z}$
in the imaginary quadatic field ${\Bbb Q}(\sqrt{-D})$; the multiplication by $\alpha\in {\goth R}$ generates an
endomorphism  $(a,b,c,d)\in M_2({\Bbb Z})$ of the lattice $L_{CM}$. 
It is known, that the endomorphisms of lattice $L_{CM}$ and endomorphisms of the
pseudo-lattice $\Lambda_{RM}=F(L_{CM})$ are related by the following explicit 
map  [4,  p. 524]:
\begin{equation}\label{eq8}
\left(\matrix{a & b\cr c & d}\right)\in End~(L_{CM})
\longmapsto
 \left(\matrix{a & b\cr -c & -d}\right)\in End~(\Lambda_{RM}),
\end{equation}
Moreover, one can always assume  $d=0$ in a proper basis of $L_{CM}$. 
We shall consider the following two cases.

\medskip
{\sf Case I.} If $D\equiv 1~mod~4$ then by (\ref{eq1})
${\goth R}={\Bbb Z}+({f+\sqrt{-f^2D}\over 2}){\Bbb Z}$;
thus the multiplier $\alpha={2m+fn\over 2}+\sqrt{{-f^2Dn^2\over 4}}$ for some 
$m, n\in {\Bbb Z}$.  Therefore  multiplication by $\alpha$ corresponds to an  endomorphism $(a,b,c,0)\in M_2({\Bbb Z})$,
where 
\begin{equation}\label{eq9}
\left\{
\begin{array}{ccc}
a &=& Tr (\alpha)=\alpha+\bar\alpha=2m+fn\\
b &=& -1\\
c &=& N (\alpha)=\alpha\bar\alpha=\left({2m+fn\over 2}\right)^2+{f^2Dn^2\over 4}. 
\end{array}
\right.
\end{equation}
To calculate a primitive generator of endomorphisms of the lattice $L_{CM}$
one should  find a multiplier $\alpha_0\ne 0$ such that
\begin{equation}
|\alpha_0|=\min_{m.n\in {\Bbb Z}}|\alpha|=\min_{m.n\in {\Bbb Z}}\sqrt{N(\alpha)}.
\end{equation}
From the last equation of (\ref{eq9}) the minimum is attained for $m=-{f\over 2}$ and
$n=1$ if $f$ is even or  $m=-f$ and $n=2$ if $f$ is odd. Thus
\begin{equation}\label{eq11}
\alpha_0=\cases{\pm{f\over 2}\sqrt{-D}, & if $f$ is even\cr
                \pm f\sqrt{-D}, & if $f$ is odd.}
\end{equation}
To find the matrix form of the endomorphism $\alpha_0$,  we shall substitute
in (\ref{eq8})  $a=d=0, b=-1$ and $c={f^2D\over 4}$ if $f$ is even or
$c=f^2D$ if $f$ is odd.  Thus the Teichm\"uller functor maps the
multiplier $\alpha_0$ into
\begin{equation}\label{eq12}
F(\alpha_0)=\cases{\pm{f\over 2}\sqrt{D}, & if $f$ is even\cr
                \pm f\sqrt{D}, & if $f$ is odd.}
\end{equation}
Comparing equations (\ref{eq11}) and (\ref{eq12}) one verifies
that formula  $F(E_{CM}^{(-D,f)})=A_{RM}^{(D,f)}$ is true in this case.

\medskip
{\sf Case II.} If $D\equiv 2$ or $3 ~mod~4$ then by (\ref{eq1}) 
${\goth R}={\Bbb Z}+(\sqrt{-f^2D})~{\Bbb Z}$;
thus the multiplier $\alpha=m+\sqrt{-f^2Dn^2}$ for some 
$m, n\in {\Bbb Z}$.  A  multiplication by $\alpha$ corresponds to an  
endomorphism $(a,b,c,0)\in M_2({\Bbb Z})$,  where 
\begin{equation}\label{eq13}
\left\{
\begin{array}{ccc}
a &=& Tr (\alpha)=\alpha+\bar\alpha=2m\\
b &=& -1\\
c &=& N (\alpha)=\alpha\bar\alpha=m^2+f^2Dn^2. 
\end{array}
\right.
\end{equation}
We shall repeat the argument of {\sf Case I}; then
from the last equation of (\ref{eq13}) the minimum of $|\alpha|$ is attained for $m=0$ and
$n=\pm 1$. Thus $\alpha_0=\pm f\sqrt{-D}$.

To find the matrix form of the endomorphism $\alpha_0$ we substitute
in (\ref{eq8})  $a=d=0, b=-1$ and $c=f^2D$.  Thus the Teichm\"uller functor maps the
multiplier $\alpha_0=\pm f\sqrt{-D}$  into $F(\alpha_0)=\pm f\sqrt{D}$.  
 In other words, formula  $F(E_{CM}^{(-D,f)})=A_{RM}^{(D,f)}$ is true in this case
as well.

\medskip
Since all possible cases are exhausted, lemma \ref{lm3} is proved.
$\square$

\begin{lem}\label{lm4}
For every $N\ge 1$  there exists a holomorphic map $X_1(N)\to X(N)$.
\end{lem}
{\it Proof.} Indeed, $\Gamma(N)$ is a normal subgroup of index $N$  of the
group $\Gamma_1(N)$;  therefore there exists a degree $N$ holomorphic 
map $X_1(N)\to X(N)$.
$\square$

\bigskip
Theorem \ref{thm1} follows from lemmas \ref{lm1}-\ref{lm3} and lemma \ref{lm4} for
$N=fD$. 
$\square$

\begin{rmk}
\textnormal{
While this note was in print,  the author came across a preprint (D'Andrea, Fiore \& Franco,  2013 \cite{DanFioFra1}).
Using the idea of quantum deformation of the line bundles over elliptic curves,  the authors 
establish a remarkable formula
\begin{equation}\label{eq15}
\tau-{p\theta\over 2}i\in {\Bbb Z}+{\Bbb Z}i,
\end{equation}
where $p\in {\Bbb Z}$ is the first Chern class of the line bundle.  
The reader is encouraged to verify,  that  theorem \ref{thm1}  satisfies 
equation (\ref{eq15}) for a  line bundle of the Chern class $p=2f'$ with 
$\tau=f\sqrt{-D}$ and $\theta=\sqrt{D}$.    
}
\end{rmk}

\bigskip\noindent
{\sf Acknowledgment.} 
I  thank the referee for helpful comments.



\vskip1cm

\textsc{The Fields Institute for Research in Mathematical Sciences, Toronto, ON, Canada,  
E-mail:} {\sf igor.v.nikolaev@gmail.com}

\smallskip
{\it Current address: 1505-657 Worcester St.,  Southbridge,  MA 01550,  U.S.A.}


\begin{thebibliography}{100}

\bibitem{BS}
Z.~I.~Borevich and I.~R.~Shafarevich, Number Theory,  Acad. Press, 1966.



\bibitem{DanFioFra1}
F.~D'Andrea, G.~Fiore  and  D.~Franco,   Modules over the noncommutative torus
and elliptic curves,    {\sf arXiv:1307.6802}



\bibitem{Hec1}
E.~Hecke, Bestimmung der Perioden gewisser Integrale durch die Theorie
der Klassenk\"orper,  Math.~Z. 28 (1928), 708-727. 


\bibitem{Man1}
Yu.~I.~Manin, Real multiplication and noncommutative geometry,
in ``Legacy of Niels Hendrik Abel'', 685-727, Springer, 2004. 



\bibitem{Nik1}
I.~Nikolaev, Remark on the rank of elliptic curves, 
Osaka J. Math. 46 (2009), 515-527. 



\bibitem{Rie1}
M.~A.~Rieffel, $C^*$-algebras associated with irrational rotations,
Pacific J. of Math. 93 (1981), 415-429.




\bibitem{S}
J.~H.~Silverman, Advanced Topics in the Arithmetic of Elliptic Curves,
GTM 151, Springer 1994.




\end{thebibliography}
\end{document}